\documentclass[preprint]{elsarticle}
\usepackage{amssymb}
\usepackage{amsfonts}
\usepackage{amsmath,amsthm}
\usepackage{graphicx}
\usepackage{graphics}
\usepackage{array}
\usepackage{subfigure}
\journal{\textbf{.}}
\oddsidemargin 0cm \numberwithin{equation}{section}

\begin{document}
\title{
A new operational matrix based on Bernoulli polynomials
}
\author[1]{J.A. Rad\corref{cor1}}
\ead{j.amanirad@gmail.com;j\_amanirad@sbu.ac.ir}
\author[2]{S. Kazem}
\author[3]{M. Shaban}
\author[1]{K. Parand}
\ead{k\_parand@sbu.ac.ir}
\cortext[cor1]{Corresponding author}
\address[1]{Department of
Computer Sciences, Faculty of Mathematical Sciences, Shahid Beheshti University, Evin, P.O. Box 198396-3113,Tehran,Iran}
\address[2]{Department of Applied Mathematics, Faculty of Mathematics and Computer Science, Amirkabir University of Technology,
No. 424, Hafez Ave., 15914, Tehran, Iran}
\address[3]{Department of Scientific Computing, Florida State University, 400 Dirac Science Library, Tallahassee, FL 32306, USA}

\begin{abstract}
In this research, the Bernoulli polynomials are introduced. The properties of these polynomials are employed to construct the operational matrices of integration together with the derivative and product.
These properties are then utilized to transform the differential equation to a matrix
equation which corresponds to a system of algebraic equations with unknown
Bernoulli coefficients. This method can be used for many problems such as differential equations, integral equations and so on.
Numerical examples show the method is computationally simple and also illustrate the efficiency and accuracy of the method.
 \end{abstract}
\begin{keyword} Bernoulli polynomial; Operational matrix; Galerkin method.\\
 AMS subject classification: 65M70, 65L60.
\end{keyword}
\maketitle
\section{Introduction}
Differential equations and their solutions play a major role in science and engineering.
A physical event can be modeled by the differential equation, an integral equation or an integro-differential equation or a system of these
equations.
Since few of these equations cannot be solved explicitly, it is often necessary to resort to the numerical techniques
which are appropriate combinations of numerical integration and interpolation. The solution of this equations occurring in physics, biology and engineering are based on numerical methods such as the Runge-Kutta methods. In recent years, the differential, integral and integro-differential equations have been solved
using the homotopy perturbation method \cite{Yildirim..Kaplan.HPM,He.wave.HPM}, the radial basis functions method \cite{Kazem.amani.RBF.EABE}, the collocation method \cite{parand.JCP}, the Homotopy analysis method \cite{Abbasbandy.Shivanian.HAM.INT,Liao.HAM}, the Tau method \cite{Parand.Razzaghi.scripta,Parand.Razzaghi.IJCM}, the Variational iteration method \cite{Yousefi.Dehghan.IJCM,wazwaz.VIM}, the Legendre matrix method \cite{Saadatmandi.Dehghan.Leg.Tau,Saadatmandi.Dehghan.Leg.Tau.2}, the Haar Wavelets operational matrix \cite{Gu.Jiang.Harr}, the Shifted Chebyshev direct method \cite{Horng.Chou}, the Legendre Wavelets operational matrix \cite{Razzaghi.yousefi.wavelet}, Sine-Cosine Wavelets operational matrix \cite{Razzaghi.yousefi.sin.cos}, the operational matrices of Bernstein polynomials \cite{S.A.Yousefi.Behroozifar.2010} and so on.\\
Polynomial series and orthogonal functions have
received considerable attention in dealing with various
problems of differential, integral and integro-differential equations.
Polynomials are incredibly useful mathematical tools as they are simply defined, can be calculated
quickly on computer systems and represent a tremendous variety of functions. They can be differentiated
and integrated easily, and can be pieced together to form spline curves that can approximate any function to
any accuracy desired.
Also, the main characteristic
of this technique is that it reduces these problems to
those of solving a system of algebraic equations, thus
greatly simplifies the problems.\\
Bernoulli polynomials play an important role in various expansions and
approximation formulas which are useful both in analytic theory of numbers
and in classical and numerical analysis.
The Bernoulli polynomials and numbers have been generalized by
Norlund \cite{Norlund} to the Bernoulli polynomials and numbers of higher order.
Also, Vandiver in \cite{Vandiver} generalized the Bernoulli numbers. Analogous
polynomials and sets of numbers have been defined from time
to time, witness the Euler polynomials and numbers and the so-called
Bernoulli polynomials of the second kind.
These polynomials can be defined by
various methods depending on the applications \cite{Cheon.Bernoulli,Burak.Simsek.Bernoul,Lu.Bernoulli,Agoh.Dilcher.Bernoulli,Buric.Elezovic.Bernoulli,Natalini.Bernoulli,Kurt.Bernoulli}.\\
In this paper, the Galerkin method \cite{Gottlieb.Hussaini} based on operational matrices of integration, differentiation
and product for the Bernoulli polynomials are presented.\\\\
The remainder of this paper is organized as follows:
In Section \ref{bernoulli},
we describe the basic formulation of Bernoulli polynomials
required for our subsequent development. Section \ref{fun.app}
is devoted to the function approximation by using
this polynomials basis. In Sections \ref{oper.matrix}, we explain
general procedure of forming of operational matrices
of integration, differentiation and product.
In Section \ref{exam}, we report our numerical findings and
demonstrate the validity, accuracy and applicability
of the operational matrices by considering numerical
examples. Also a conclusion is given in
the last Section.
\section{Bernoulli polynomials}\label{bernoulli}
The classical Bernoulli polynomials of $n$-th degree are defined on the interval [0, 1] as
\begin{eqnarray}\nonumber
B_{n}(x)=\sum_{k=0}^{n}{n \choose k}B_{k}x^{n-k}~,
\end{eqnarray}
where $B_{k} := B_{k}(0)$ is the Bernoulli number for each $k = 0 , 1 , . . . , n$.
Thus, the first four such polynomial, respectively, are
\begin{eqnarray}\nonumber
&&B_{0}(x)=1~,~~~~~~~~~~~~~~~~B_{1}(x)=x-\frac{1}{2}~,\\\nonumber
&&B_{2}(x)=x^2-x+\frac{1}{6}~,~~~B_{3}(x)=x^3-\frac{3}{2}x^2+\frac{1}{2}x~,
\end{eqnarray}
Leopold Kronecker expressed the Bernoulli number $B_{n}$ in the following form
\begin{eqnarray}\nonumber
B_{n}=-\sum_{k=1}^{n+1}\frac{(-1)^k}{k}{n+1 \choose k}\sum_{j=1}^{k}j^{n}~,
\end{eqnarray}
for $n\geq0~,~~n\neq1$. If $n=1$, we defined $B_{1}=-\frac{1}{2}$.\\
Determinant form of the Bernoulli polynomial of $n$-th degree is defined by F. Costabile \cite{Costabile}
\begin{eqnarray}\nonumber
B_{0}(x)=1~,
\end{eqnarray}
and
\begin{align}\nonumber
B_{n}(x)=\frac{(-1)^{n}}{(n-1)!}
&\left|
 \begin{array}{ccccccc}
1 & x & x^2 & x^3 & \dots  & x^{n-1} & x^n\\
1 & \frac{1}{2} & \frac{1}{3} & \frac{1}{4} & \dots  & \frac{1}{n} & \frac{1}{n+1}\\
0 & 1 & 1 & 1 & \dots  & 1 & 1\\
0 & 0 & 2 & 3 & \dots  & n-1 & n\\
0 & 0 & 0 & {3 \choose 2} & \dots  & {n-1 \choose 2} & {n \choose 2}\\
\vdots          & \vdots  & \vdots & \vdots         & \ddots & \vdots   & \vdots       \\
0 & 0 & 0 & 0 & \dots  & {n-1 \choose n-2} & {n \choose n-2}\\
\end{array}
\right|~.
\end{align}
for each $n = 1, 2, ...~$.\\
The Bernoulli polynomials play an important role in different areas of mathematics, including number theory and the theory of finite differences.
This polynomial have many similar properties. These polynomials produce the following exponential generating function
\begin{eqnarray}\label{EGF}
\frac{te^{xt}}{e^{t}-1}=\sum_{n=0}^{\infty}B_{n}(x)\frac{t^{n}}{n!}~,
\end{eqnarray}
Since they satisfy the well-known relation
\begin{eqnarray}\label{der.Ber}
\frac{d}{dx}B_{n}(x)=nB_{n-1}(x)~,
\end{eqnarray}
(for all $n \geq 1$), which follows easily from Eq. (\ref{EGF}), it is to be expected that integrals figure prominently in the study of these
polynomials. The most immediate integral formula is obtained by integrating Eq. (\ref{der.Ber})
\begin{eqnarray}\label{Int.Ber}
B_{n}(x)=n~\int_{0}^{x}B_{n-1}(t)~dt+B_{n}~,
\end{eqnarray}
One the other of these properties is
\begin{eqnarray}\label{11}
&&B_{n}(x+1)-B_{n}(x)=nx^{n-1}~,\\\label{22}
&&B_{n}(x+1)=\sum_{k=0}^{n}{n \choose k}B_{k}(x)~,\\\nonumber
&&B_{n}(1-x)=(-1)^{n}B_{n}(x)~,\\\nonumber
&&(-1)^{n}B_{n}(-x)=B_{n}(x)+nx^{n-1}~,\\\nonumber
&&\int_{0}^{1}B_{n}(x)~dx=0~,~~~n\geq1~.
\end{eqnarray}
From Eq. (\ref{11}) and (\ref{22}), we obtain for any integer $n \geq 0$,
\begin{eqnarray}\label{inver}
\sum_{k=0}^{n}{n+1 \choose k}B_{k}(x)=(n+1)x^{n}~.
\end{eqnarray}
Also the following interesting integral for a product of two
Bernoulli polynomials appears in the book by N\"{o}rlund for all $k +m \geq 2$,
\begin{eqnarray}\label{prooo}
\int_{0}^{1}B_{n}(t)B_{m}(t)~dt=(-1)^{n-1}\frac{m!n!}{(m+n)!}B_{n+m}~.
\end{eqnarray}
It can be easily shown that any given polynomial of degree $n$ can be expanded in terms
of linear combination of the basis functions
\begin{eqnarray}\label{ber}
P(x)=\sum_{i=0}^{n}C_{i}B_{i}(x)=C^T\mathbf{B}(x)~,
\end{eqnarray}
where $C$ and $\mathbf{B}(x)$ are $(n+1)\times 1$ vectors given by
\begin{eqnarray}\nonumber
&&C=[C_0, C_1, . . . , C_
{n}]^T,\\\label{vector.B.Ber}
&&\mathbf{B}(x)=[B_{0}(x),B_{1}(x)
, . . . , B_{n}(x)]^T~.
\end{eqnarray}
For each $i=0,1,...,n~$, we can obtain the following matrix form of $B_{i}(x)$
\begin{eqnarray}\nonumber
&&B_{i}(x)=\sum_{k=0}^{i}{i \choose k}B_{k}x^{i-k}\\\nonumber
&&={i \choose i}B_{i}+{i \choose i-1}B_{i-1}x+...+{i \choose 1}B_{1}x^{i-1}+{i \choose 0}B_{0}x^{i}\\\nonumber
&&=
\left[
 \begin{array}{ccccccccccc}
{i \choose i}B_{i} & {i \choose i-1}B_{i-1} & {i \choose i-2}B_{i-2} & \dots  & {i \choose 1}B_{1} & {i \choose 0}B_{0} & \overbrace{0 ~~ 0 ~~ \dots ~~ 0 ~~ 0}^{n-i}\\
\end{array}
\right]~
\left[
 \begin{array}{c}
1\\
x\\
x^2\\
\vdots\\
x^i\\
x^{i+1}\\
\vdots\\
x^{n}
\end{array}
\right]~\\
&&=\mathbf{M}_{i}~\mathrm{T}(x)~.
\end{eqnarray}
where
\begin{eqnarray}\nonumber
\mathrm{T}(x)=\Bigg[1~~~x~~~x^2~~~~...~~~~x^n\Bigg]^{T}~
\end{eqnarray}
and
\begin{eqnarray}\nonumber
\mathbf{M}_{i}=
\left[
 \begin{array}{ccccccccccc}
{i \choose i}B_{i} & {i \choose i-1}B_{i-1} & {i \choose i-2}B_{i-2} & \dots  & {i \choose 1}B_{1} & {i \choose 0}B_{0} & \overbrace{0 ~~ 0 ~~ \dots ~~ 0 ~~ 0}^{n-i}\\
\end{array}
\right]^{T}~
\end{eqnarray}
Now, we can expand the matrix $\mathbf{B}(x)=[B_{0}(x),B_{1}(x)
, . . . , B_{n}(x)]^T$ as
\begin{eqnarray}\nonumber
&&\mathbf{B}(x)=[B_{0}(x),B_{1}(x)
, . . . , B_{n}(x)]^T\\\nonumber
&&~~~=[\mathbf{M}_{0}~\mathrm{T}(x),\mathbf{M}_{1}~\mathrm{T}(x), . . . , \mathbf{M}_{n}~\mathrm{T}(x)]^T\\\nonumber
&&~~~=[\mathbf{M}_{0},\mathbf{M}_{1}, . . . , \mathbf{M}_{n}]^T~\mathrm{T}(x)\\\nonumber
&&~~~=
\left[
 \begin{array}{ccccccc}\label{matrice}
B_{0} & 0 & 0 & 0 & 0 & \dots  & 0 \\
{1 \choose 1}B_{1} & {1 \choose 0}B_{0} & 0 & 0 & 0 & \dots  & 0 \\
{2 \choose 2}B_{2} & {2 \choose 1}B_{1} & {2 \choose 0}B_{0} & 0 & 0 & \dots  & 0 \\
{3 \choose 3}B_{3} & {3 \choose 2}B_{2} & {3 \choose 1}B_{1} & {3 \choose 0}B_{0} & 0 & \dots  & 0 \\
\vdots     & \vdots     & \vdots   & \vdots & \vdots       & \ddots & \vdots \\
{n \choose n}B_{n} & {n \choose n-1}B_{n-1} & {n \choose n-2}B_{n-2} & {n \choose n-3}B_{n-3} & 0 & \dots  &  {n \choose 0}B_{0}\\
\end{array}
\right]~\mathrm{T}(x)\\\nonumber\\
&&~~~=\mathbf{M}~\mathrm{T}(x)
\end{eqnarray}
Then, one has
\begin{eqnarray}\label{B.to.T}
\mathbf{B}(x)=\mathbf{M}~\mathrm{T}(x)~,
\end{eqnarray}
where $\mathbf{M}$ is the $(n+1)\times(n+1)$ matrix and
\begin{eqnarray}\nonumber
\{\mathbf{M}\}_{i,j=1}^{n+1}
=
\begin{cases}
B_{i-j}{i-1 \choose i-j}~, & i\geq j~,\\
0~, & i < j~.
\end{cases}
\end{eqnarray}
Matrix $\mathbf{M}$ is an lower triangular matrix and $\det(\mathbf{M})=1$, so $\mathbf{M}$ is an invertible matrix.
On the other hand by using Eq.(\ref{inver}) we have
\begin{eqnarray}\nonumber
&&\sum_{k=0}^{n}{n+1 \choose k}B_{k}(x)=(n+1)x^{n}~,\\\nonumber\\\label{T.to.B}
&&\mathbf{Q}\mathbf{B}(x)=\mathrm{T}(x)~~~~~\Rightarrow ~~~~~\mathbf{B}(x)=\mathbf{Q}^{-1}~\mathrm{T}(x)~,
\end{eqnarray}
where $\mathbf{Q}$ is the $(n+1)\times(n+1)$ matrix and
\begin{eqnarray}\nonumber
\{\mathbf{Q}\}_{i,j=1}^{n+1}
=
\begin{cases}
\frac{1}{i}{i \choose j-1}~, & i\geq j~,\\
0~, & i < j~.
\end{cases}
\end{eqnarray}
Matrix $\mathbf{Q}$ is an lower triangular matrix and $\det(\mathbf{Q})=1$, so $\mathbf{Q}$ is an invertible matrix.
Then  from Eqs. (\ref{B.to.T}) and (\ref{T.to.B}), one has
\begin{eqnarray}\nonumber
\mathbf{Q}^{-1}=\mathbf{M}~,
\end{eqnarray}
\section{Function approximation}\label{fun.app}
Let us define $\Lambda=\{x|~~0 \leq x \leq 1\}$ and\\\\
$~~~~~L^{2}(\Lambda)=\{\nu:\Lambda\rightarrow\mathbb{R}|\nu$ is measurable and $||\nu||_2<\infty\}~,$\\\\
where
\begin{eqnarray}\label{normm}
||\nu||_2=\Bigg(\int_{0}^{1}|\nu(x)|^2~dx\Bigg)^{1/2}~,
\end{eqnarray}
is the norm induced by the inner product of the space $L^2(\Lambda)$,
\begin{eqnarray}\label{product}
<u,\nu>=\int_{0}^{1}u(x)\nu(x)~dx~,
\end{eqnarray}
Now, we suppose
\begin{eqnarray}\nonumber
\mathfrak{B}_{N}=\{B_{0}(x),B_{1}(x),...,B_{N}(x)\}~,
\end{eqnarray}
where $\mathfrak{B}_{N}$ is finite dimensional subspace, therefore $\mathfrak{B}_{N}$ is a complete subspace of $L^2(\Lambda)$.\\
The interpolating function of a smooth function $u$ on a finite interval is denoted by $\xi_N u$. It is an element of $\mathfrak{B}_{N}$ and
\begin{eqnarray}\nonumber
\xi_N u=\sum_{k=0}^{N}a_kB_{k}(x)~,
\end{eqnarray}
$\xi_N u$ is the best projection of $u$ upon $\mathfrak{B}_{N}$ with respect to the inner product Eq. (\ref{product}) and the norm Eq. (\ref{normm}). Thus, we have
\begin{eqnarray}\nonumber
<\xi_N u-u,B_{i}(x)>=0~,~~~~~~~\forall~~B_{i}(x) \in \mathfrak{B}_{N}~,
\end{eqnarray}
 or equivalently
\begin{eqnarray}\label{func.appr}
\xi_N u(x)=\sum_{i=0}^{N}a_iB_{i}(x)=A^{T}\mathbf{B}(x)~,
\end{eqnarray}
where $A$ and $\mathbf{B}(x)$ are $(N+1)\times 1$ vectors given by
\begin{eqnarray}\nonumber
&&A=[a_0, a_1, . . . , a_
{N}]^T,\\\label{vector.B.Ber}
&&\mathbf{B}(x)=[B_{0}(x),B_{1}(x)
, . . . , B_{N}(x)]^T~,
\end{eqnarray}
therefore $A$ can be obtained by
\begin{eqnarray}\nonumber
A^{T}=<\xi_N u(x),\mathbf{B}^{T}(x)><\mathbf{B}(x),\mathbf{B}^{T}(x)>^{-1}~,
\end{eqnarray}
where
\begin{eqnarray}\nonumber
<\xi_N u(x),\mathbf{B}^{T}(x)>=\int_{0}^{1}\xi_N u(x)~\mathbf{B}^{T}(x)~dx~,
\end{eqnarray}
is an $1 \times (N+1)$ vector and
\begin{eqnarray}\label{Galer.D}
&&\mathbf{D}=<\mathbf{B}(x),\mathbf{B}^{T}(x)>=\int_{0}^{1}\mathbf{B}(x)~\mathbf{B}^{T}(x)~dx~,
\end{eqnarray}
is an $(N+ 1) \times (N+ 1)$ matrix and is said
dual matrix of $\mathbf{B}(x)$. Suppose that $\mathbf{D}=\{d_{i,j}\}_{i,j=1}^{N+1}$, and by using Eq. (\ref{prooo}) one has
\begin{eqnarray}\nonumber
&&d_{i,j}=\int_{0}^{1}B_{i}(x)B_{j}(x)~dx=(-1)^{i-1}\frac{i!j!}{(i+j)!}B_{i+j}~,
\end{eqnarray}
which this shows the matrix $D$ is symmetric and invertible.
For example, if $N=5$, we have
\begin{eqnarray}\nonumber
\mathbf{D}=
\left[
 \begin{array}{cccccc}\label{matrice}
1 & 0 & 0 & 0 & 0 & 0 \\
0 & \frac{1}{12} & 0 & \frac{-1}{120} & 0 & \frac{1}{252} \\
0 & 0 & \frac{1}{180} & 0 & \frac{-1}{630} & 0 \\
0 & \frac{-1}{120} & 0 & \frac{1}{840} & 0 & \frac{-1}{1680} \\
0 & 0 & \frac{-1}{630} & 0 & \frac{1}{2100} & 0 \\
0 & \frac{1}{252} & 0 & \frac{-1}{1680} & 0 & \frac{5}{16632} \\
\end{array}
\right]
\end{eqnarray}
\section{The operational matrix}\label{oper.matrix}
It is known that operational matrices are employed for solving many engineering and physical problems such as
dynamical systems [5], optimal control systems [6–8], robotic systems [9] etc. Furthermore they are used in several areas
of numerical analysis, and they hold particular importance in various subjects such as integral equations [10], differential
equations [11,12], calculus of variations [13], partial differential equations [14], integro-differential equations [15,16] etc.
Also many books and papers have employed the operational matrix for spectral methods [17,18]. \\
Let us start this section
by introducing operational matrices. Suppose that
\begin{eqnarray}\nonumber
\mathbf{B}(x)=[B_{0}(x),B_{1}(x)
, . . . , B_{N}(x)]^T~,
\end{eqnarray}
the matrices $\mathcal{D}$ and $\mathcal{I}$ are named respectively as the operational matrices of derivatives and integrals if and only if
\begin{eqnarray}\nonumber
&&\frac{d}{dx}\mathbf{B}(x)=\mathcal{D}~\mathbf{B}(x)~,\\\label{oper.integral}
&&\int_{0}^{x}\mathbf{B}(t)~dt\simeq\mathcal{I}~\mathbf{B}(x)~.
\end{eqnarray}
Furthermore assume $C = [c_0, c_1, ... , c_N]$, $\tilde{C}$ is named as the operational matrix of the product if and only if
\begin{eqnarray}\nonumber
\mathbf{B}(x)\mathbf{B}^{T}(x)C\simeq\tilde{C}~\mathbf{B}(x)~.
\end{eqnarray}
\subsection{Operational matrix of derivative}
\textbf{Theorem 1}. Let $\mathbf{B}(x)$ be the Bernoulli vector, $\mathcal{D}$ is the $(N+1) \times (N+1)$ operational matrix of derivatives,
then the elements of $\mathcal{D}$ are obtained as
\begin{eqnarray}\nonumber
\{D_{i,j}\}_{i,j=0}^{N}=
\begin{cases}
i~, & i=j+1~, \\
0~, & otherwise~.
\end{cases}
\end{eqnarray}
\textbf{Proof}.
By using Eq. (\ref{der.Ber}), we have
\begin{eqnarray}\nonumber
&&\frac{d}{dx}\mathbf{B}(x)=\Bigg[\frac{d}{dx}B_{0}(x),\frac{d}{dx}B_{1}(x),...,\frac{d}{dx}B_{N}(x)\Bigg]\\\nonumber
&&~~~~~=\Bigg[0,B_{0}(x),2B_{1}(x),3B_{2}(x),...,NB_{N-1}(x)\Bigg]\\\nonumber
&&~~~~~=
\left[
 \begin{array}{ccccccc}\label{matrice}
0 & 0 & 0 & 0 & \dots & 0 & 0\\
1 & 0 & 0 & 0 & \dots & 0 & 0\\
0 & 2 & 0 & 0 & \dots & 0 & 0\\
0 & 0 & 3 & 0 & \dots & 0 & 0\\
\vdots     & \vdots     & \vdots   & \vdots   & \ddots & \vdots &\vdots\\
0 & 0 & 0 & 0 & \dots & N & 0\\
\end{array}
\right]
~\mathbf{B}(x)=\mathcal{D}~\mathbf{B}(x),
\end{eqnarray}

\subsection{Operational matrix of integration}
Let $\mathbf{B}(x)$ be the Bernoulli vector, $\mathcal{I}$ is the $(N+1) \times (N+1)$ operational matrix of integrative,
then by using Eq. (\ref{Int.Ber})
\begin{eqnarray}\nonumber
\int_{0}^{x}\mathbf{B}(t)~dt\simeq\mathcal{I}~\mathbf{B}(x)~.
\end{eqnarray}

\begin{eqnarray}\nonumber
\int_{0}^{x}\mathbf{B}(t)~dt=
\left[
 \begin{array}{c}\label{matrice}
\int_{0}^{x}B_{0}(t)~dt\\
\int_{0}^{x}B_{1}(t)~dt\\
\vdots\\
\int_{0}^{x}B_{N}(t)~dt \\
\end{array}
\right]
~=
\left[
 \begin{array}{c}\label{matriceee}
B_{1}(x)-B_{1}\\
\frac{B_{2}(x)-B_{2}}{2}\\
\frac{B_{3}(x)-B_{3}}{3}\\
\vdots\\
\frac{B_{N+1}(x)-B_{N+1}}{N+1}\\
\end{array}
\right]~,
\end{eqnarray}
by using expansion of $B_{i}(x)$, we have
\begin{eqnarray}\nonumber
&&\frac{B_{i}(x)-B_{i}}{i}=\frac{1}{i}\sum_{k=0}^{i}{i \choose k}B_{k}x^{i-k}-\frac{1}{i}B_{i}\\\nonumber
&&=\frac{1}{i}{i \choose i-1}B_{i-1}x+\frac{1}{i}{i \choose i-2}B_{i-2}x^2+...+\frac{1}{i}{i \choose 1}B_{1}x^{i-1}+\frac{1}{i}B_{0}x^{i}\\\nonumber
&&=
\left[
 \begin{array}{cccccccccc}\label{matrice}
0 & \frac{1}{i}{i \choose i-1}B_{i-1} & \frac{1}{i}{i \choose i-2}B_{i-2} & \dots & \frac{1}{i}{i \choose 1}B_{1} & \frac{1}{i}B_{0} & \underbrace{0 ~~~ 0 ~~~ \dots ~~~ 0}_{N-i}\\
\end{array}
\right]~\mathrm{T}(x)~\\\label{app.n}
&&= U_{i}~\mathrm{T}(x)=U_{i}\mathbf{Q}\mathbf{B}(x),
\end{eqnarray}
where $i=1,2,...,N$. So, we just need to approximate $\frac{B_{N+1}(x)-B_{N+1}}{N+1}$. By using Eq. (\ref{func.appr}), we have
\begin{eqnarray}\label{app.n1}
\frac{B_{N+1}(x)-B_{N+1}}{N+1}\simeq\sum_{j=0}^{N}\chi_{j}B_{j}(x)=\Xi^{T}~\mathbf{B}(x)
\end{eqnarray}
where
\begin{eqnarray}\nonumber
\Xi~=~[\chi_0,\chi_1,...,\chi_N]^{T}~,
\end{eqnarray}
therefore, by substituting Eqs. (\ref{app.n}) and (\ref{app.n1}) in Eq. (\ref{matriceee}), we have
\begin{eqnarray}\nonumber
&&\int_{0}^{x}\mathbf{B}(t)~dt=
\left[
 \begin{array}{c}\label{matrice}
U_{1}\mathbf{Q}\mathbf{B}(x)\\
U_{2}\mathbf{Q}\mathbf{B}(x)\\
\vdots\\
U_{N}\mathbf{Q}\mathbf{B}(x)\\
\Xi^{T}~\mathbf{M}~\mathbf{Q}~\mathbf{B}(x)
\end{array}
\right]
=
\left[
 \begin{array}{c}\label{matrice}
U_{1}\\
U_{2}\\
\vdots\\
U_{N}\\
\Xi^{T}~\mathbf{M}\\
\end{array}
\right]
~\mathbf{Q}~\mathbf{B}(x)\\\nonumber\\
&&~~~~~~~~~~~~~~=\underbrace{\mathbf{U}~\mathbf{Q}}_{\mathcal{I}}~\mathbf{B}(x)
\end{eqnarray}
\begin{eqnarray}\nonumber
~~~~~~~~~~\Rightarrow~~~~~~~~\mathcal{I}=\mathbf{U}~\mathbf{Q}
\end{eqnarray}
where
\begin{eqnarray}\nonumber
\mathbf{U}=\Bigg[U_{1}~~~~U_{2}~~~~...~~~~U_{N}~~~~\Xi^{T}\mathbf{M}\Bigg]^{T}~,
\end{eqnarray}

\subsection{Operational matrix of product}
The following property of the product of two B-Polynomials vectors will also be applied
\begin{eqnarray}\label{Eq.tilde}
\mathbf{B}(x)\mathbf{B}^{T}(x)C\simeq\tilde{\mathbf{C}}\mathbf{B}(x)~,
\end{eqnarray}
where $\tilde{\mathbf{C}}$ is an $(N+1)\times(N+1)$ product operational matrix for the vector $C$.\\
Using the transformation Matrices between $\mathbf{B}(x)$ and $\mathrm{T}(x)$ polynomials of the previous section, we have
\begin{eqnarray}\nonumber
\mathbf{B}(x)\mathbf{B}^{T}(x)C=\mathbf{M}~\mathrm{T}(x)\mathrm{T}^{T}(x)~\underbrace{\mathbf{M}^{T}C}_{\mathbf{N}}
\end{eqnarray}
Now, we suppose that
\begin{eqnarray}\label{tx.telda}
\mathrm{T}(x)\mathrm{T}^{T}(x)~\mathbf{N}=\tilde{\mathbf{N}}~\mathrm{T}(x)~,
\end{eqnarray}
Then, we have
\begin{eqnarray}\nonumber
\mathbf{B}(x)\mathbf{B}^{T}(x)C=\mathbf{M}\tilde{\mathbf{N}}~\mathrm{T}(x)=\underbrace{\mathbf{M}\tilde{\mathbf{N}}~\mathbf{Q}}_{\tilde{\mathbf{C}}}~\mathbf{B}(x)~,
\end{eqnarray}
\begin{eqnarray}\nonumber
~~~~~\Rightarrow~~~~~~~\tilde{\mathbf{C}}=\mathbf{M}\tilde{\mathbf{N}}~\mathbf{Q}~,
\end{eqnarray}
now in order to achieve $\mathbf{\tilde{C}}$, it is sufficient to obtain the $\tilde{\mathbf{N}}$. By using Eq. (\ref{tx.telda}), we have
\begin{eqnarray}\nonumber
&&\mathrm{T}(x)\mathrm{T}^{T}(x)~\mathbf{N}=\left[
 \begin{array}{ccccc}\label{matrice}
1 & x & x^2 & \dots & x^N\\
x & x^2 & x^3 & \dots & x^{N+1}\\
x^2 & x^3 & x^4 & \dots & x^{N+2}\\
\vdots     & \vdots     & \vdots     & \ddots & \vdots\\
x^{N} & x^{N+1} & x^{N+2} & \dots & x^{2N}\\
\end{array}
\right]
\left[
 \begin{array}{c}\label{matrice}
n_0\\
n_1\\
n_2\\
\vdots\\
n_{N}\\
\end{array}
\right]\\\nonumber
&&=
\left[
 \begin{array}{ccccc}\label{matrice}
n_0 & n_1 & n_2 & \dots & n_N\\
0 & n_0 & n_1 & \dots & n_{N-1}\\
0 & 0 & n_0 & \dots & n_{N-2}\\
\vdots     & \vdots     & \vdots     & \ddots & \vdots\\
0 & 0 & 0 & \dots & n_0\\
\end{array}
\right]
\left[
\begin{array}{c}\label{matrice}
1\\
x\\
x^2\\
\vdots\\
x^{N}\\
\end{array}
\right]=\tilde{\mathbf{N}}~\mathrm{T}(x)~.
\end{eqnarray}

The constructed operational matrices are now used to solve the
following examples.

\section{Illustrative examples}\label{exam}
\begin{figure}
\center
\includegraphics[scale=0.5]{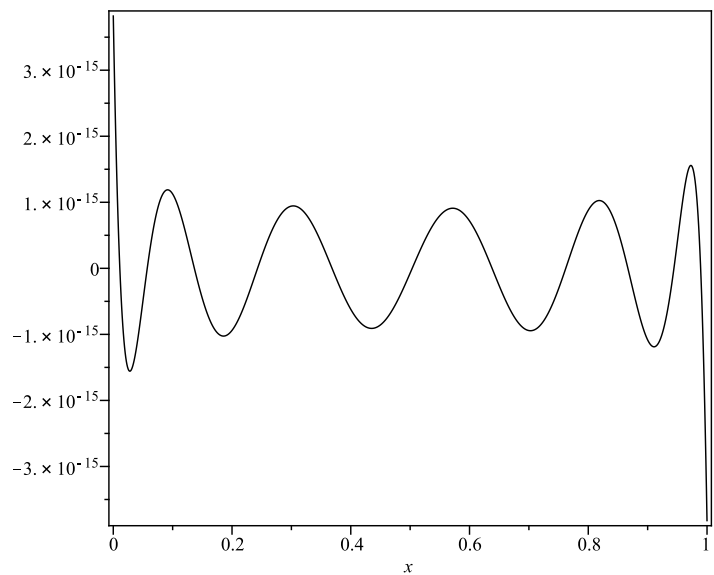}
\caption{Absolute error function for $N=10$ in Example 1.} \label{Ex1.error}
\end{figure}

\begin{figure}
\center
\includegraphics[scale=0.5]{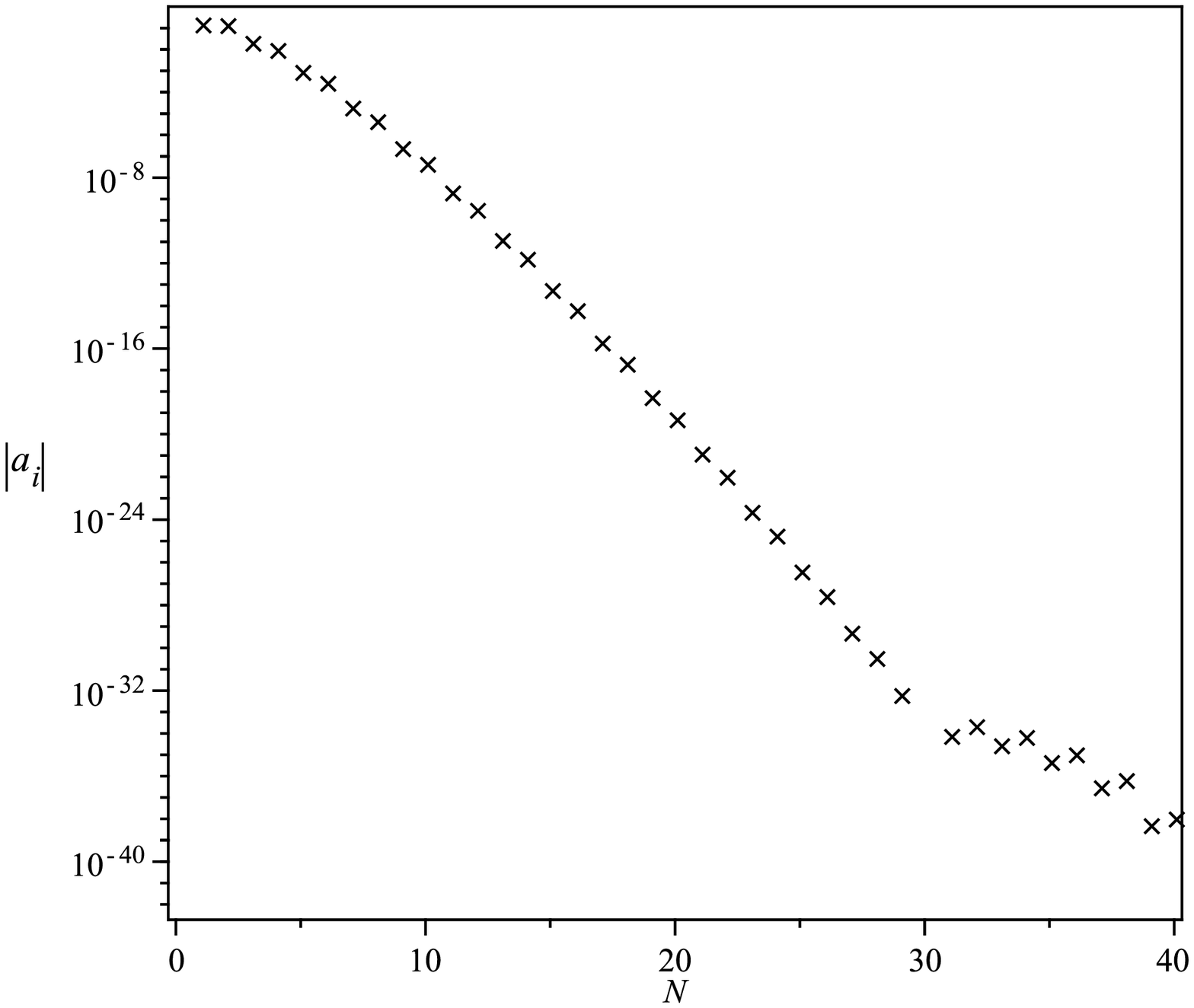}
\caption{Absolute values $|a_{i}|$ of the coefficients of the Bernoulli polynomials for some various $N$ in Example 1.} \label{Ex1.ai}
\end{figure}

\begin{figure}
\center
\includegraphics[scale=0.5]{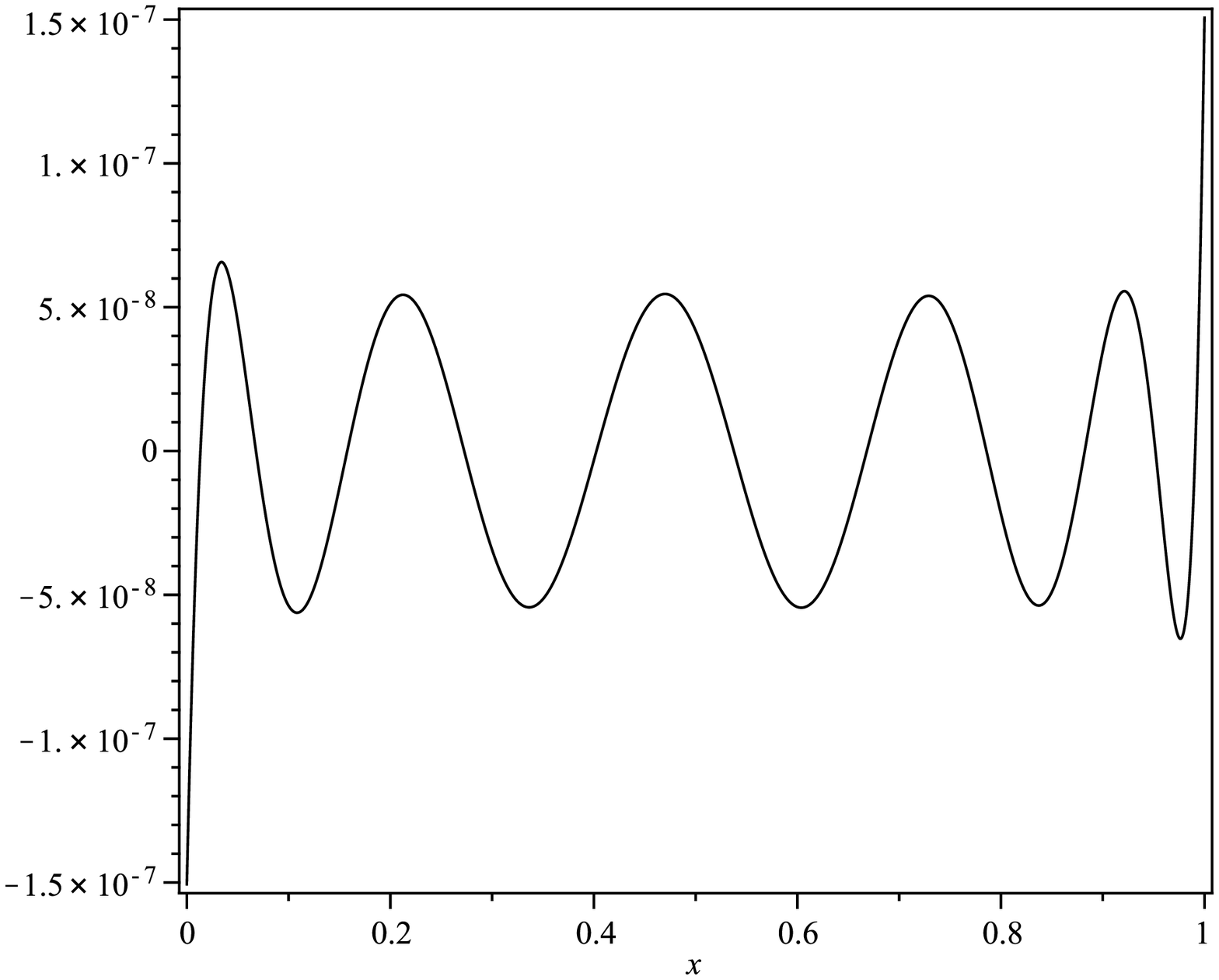}
\caption{Absolute error function for $N=10$ in Example 3.} \label{Ex3.error}
\end{figure}

\begin{figure}
\center
\includegraphics[scale=0.5]{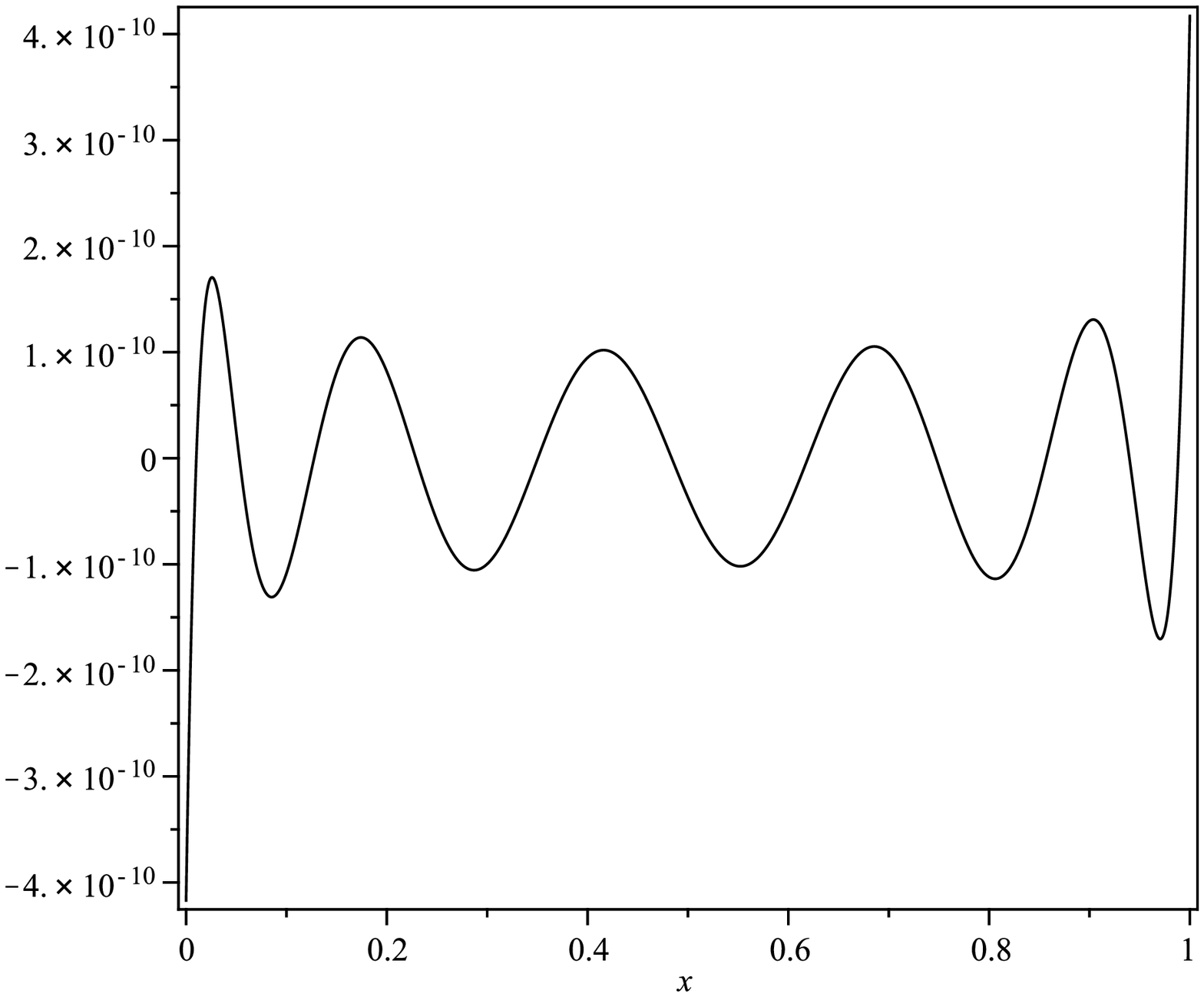}
\caption{Absolute error function for $N=10$ in Example 4.} \label{Ex4.error}
\end{figure}

\begin{figure}
\center
\includegraphics[scale=0.5]{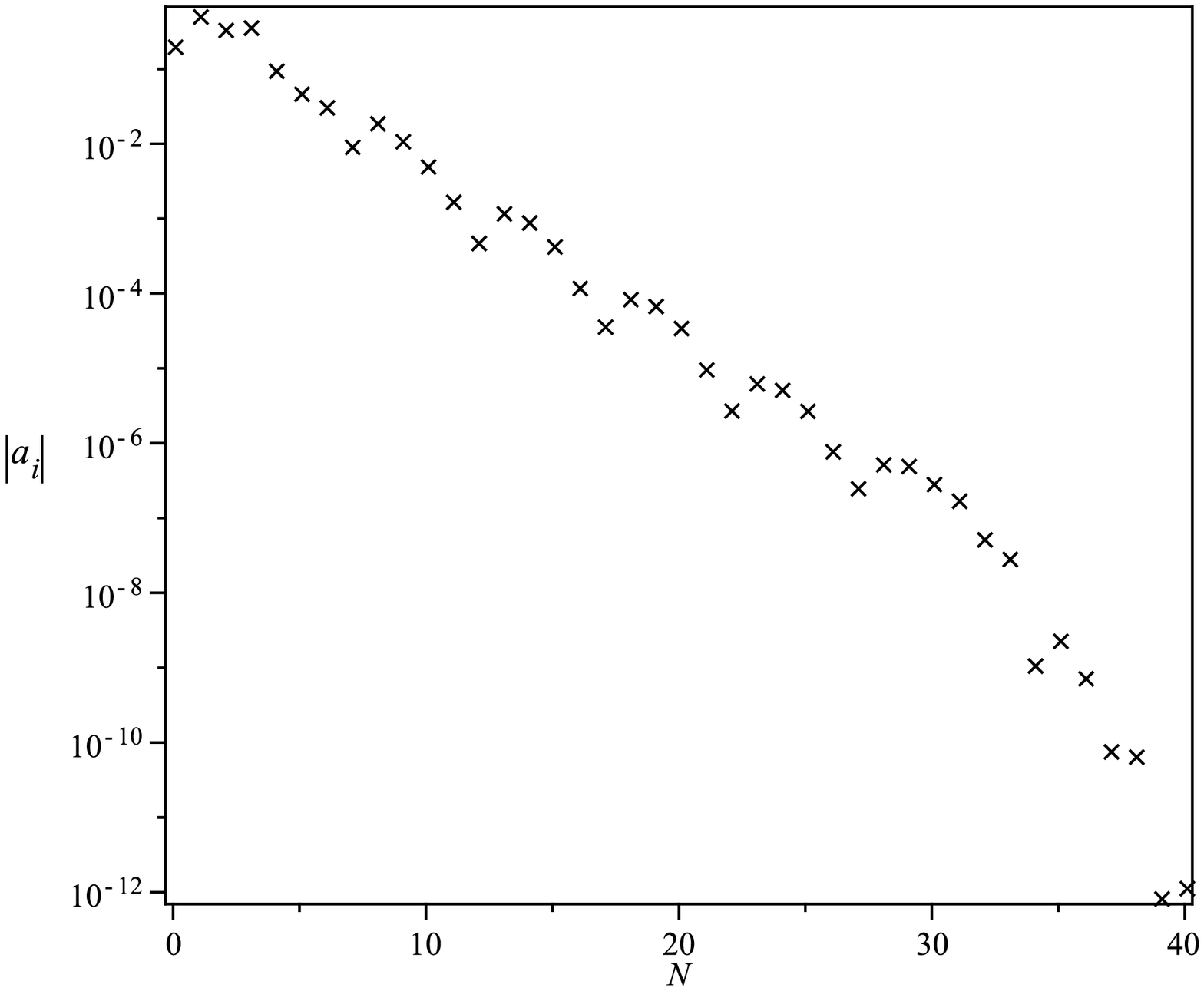}
\caption{Absolute values $|a_{i}|$ of the coefficients of the Bernoulli polynomials for some various $N$ in Example 4.} \label{Ex4.ai}
\end{figure}

 To illustrate the efficiency of the proposed method in the present paper, several test examples are carried out. The computational results obtained by using this scheme are in excellent agreement with the exact solutions. For these
comparisons the root mean square (RMS) error is applied of
the form
\begin{eqnarray}\nonumber
RMS=\sqrt{\frac{\sum_{k=1}^{M}\big(y(x_k)-y_n(x_k)\big)^2}
{M}}~,
\end{eqnarray}
where $y(x_k)$ and $y_m(x_k)$ are achieved by exact and
numerical solution on $x_k$ and $M$ is number of test points.\\\\
\textbf{Example 1.}
We consider the following Bessel differential equation
of order zero \cite{ONeil,Yousefi.Behroozifar.IJSS}
\begin{eqnarray}\label{example1.bessel}
x~\frac{d^{2}u(x)}{dx^2}+\frac{du(x)}{dx}+x~u(x)=0~,
\end{eqnarray}
with initial conditions
\begin{eqnarray}\label{initial.BEssel}
u(0)=1~,~~~~\frac{du(x)}{dx}\bigg|_{x=0}=0~,
\end{eqnarray}
The exact solution of this problem is
\begin{eqnarray}\nonumber
u(x)=\sum_{i=0}^{\infty}\frac{(-1)^{i}}{4^i(i!)^2}x^{2i}~.
\end{eqnarray}
To solve this example, we approximate $d^2u(x)/dx^2$ by the Bernoulli polynomials as
\begin{eqnarray}\nonumber
\frac{d^2u(x)}{dx^2}=A^{T}\mathbf{B}(x)~,
\end{eqnarray}
Also, by using the initial conditions Eq. (\ref{initial.BEssel}) and the operation matrix of integration Eq. (\ref{oper.integral}), we have:
\begin{eqnarray}\nonumber
&&\frac{du(x)}{dx}=A^{T}\mathcal{I}\mathbf{B}(x)~,\\\label{exe1.y}
&&u(x)=A^{T}\mathcal{I}^{2}\mathbf{B}(x)+V^{T}\mathbf{B}(x)~,
\end{eqnarray}
where $V^{T}\mathbf{B}(x)=1$ and $V=[1,0,0,...,0]^{T}$.
We can express function $x$ as
\begin{eqnarray}\nonumber
x=E^{T}\mathbf{B}(x)~,
\end{eqnarray}
where $E=[1/2,1,0,...,0]^{T}$. Therefore, equation (\ref{example1.bessel}) can be rewrite as:
\begin{eqnarray}\nonumber
E^{T}\mathbf{B}(x)\mathbf{B}^{T}(x)A+A^{T}\mathcal{I}\mathbf{B}(x)+E^{T}\mathbf{B}(x)\mathbf{B}^{T}(x)(\mathcal{I}^{2})^{T}A+E^{T}\mathbf{B}(x)\mathbf{B}^{T}(x)V=0~,
\end{eqnarray}
or
\begin{eqnarray}\nonumber
\Bigg(E^{T}\tilde{\mathbf{A}}+A^{T}\mathcal{I}+E^{T}\tilde{\mathbf{Z}}+E^{T}\tilde{\mathbf{V}}\Bigg)\mathbf{B}(x)=0~,
\end{eqnarray}
where
\begin{eqnarray}\nonumber
&&\mathbf{B}(x)\mathbf{B}^{T}(x)A=\tilde{\mathbf{A}}\mathbf{B}(x)~,\\\nonumber
&&Z=(\mathcal{I}^{2})^{T}A~,\\\nonumber
&&\mathbf{B}(x)\mathbf{B}^{T}(x)Z=\tilde{\mathbf{Z}}\mathbf{B}(x)~,\\\nonumber
&&\mathbf{B}(x)\mathbf{B}^{T}(x)V=\tilde{\mathbf{V}}\mathbf{B}(x)~,
\end{eqnarray}
Now, by using Galerkin method \cite{Gottlieb.Hussaini} and Eq. (\ref{Galer.D}), we have
\begin{eqnarray}\nonumber
&&\int_{0}^{1}\Bigg(E^{T}\tilde{\mathbf{A}}+A^{T}\mathcal{I}+E^{T}\tilde{\mathbf{Z}}+E^{T}\tilde{\mathbf{V}}\Bigg)\mathbf{B}(x)\mathbf{B}(x)~dx\\\nonumber
&&=\Bigg(E^{T}\tilde{\mathbf{A}}+A^{T}\mathcal{I}+E^{T}\tilde{\mathbf{Z}}+E^{T}\tilde{\mathbf{V}}\Bigg)\int_{0}^{1}\mathbf{B}(x)\mathbf{B}(x)~dx\\\nonumber
&&=\Bigg(E^{T}\tilde{\mathbf{A}}+A^{T}\mathcal{I}+E^{T}\tilde{\mathbf{Z}}+E^{T}\tilde{\mathbf{V}}\Bigg)\mathbf{D}=0~.
\end{eqnarray}
Since $\mathbf{D}$ is an invertible matrix, thus one has
\begin{eqnarray}\nonumber
E^{T}\tilde{\mathbf{A}}+A^{T}\mathcal{I}+E^{T}\tilde{\mathbf{Z}}+E^{T}\tilde{\mathbf{V}}=0~.
\end{eqnarray}
We generate $N + 1$ equations, therefore by solving above equations the unknown vector $A$ is achieved and
from Eq. (\ref{exe1.y}), $u(x)$ can be calculated.\\
In Table \ref{Tab1}, a comparison is made between the
approximate values using the present approach
together with the exact solution for some various $N$.
Also, the RMS errors for some various $N$ are shown in Table \ref{Tab1}.
From Table \ref{Tab1}, it can be seen that
a few term of Bernoulli polynomials is sufficient to achieve a good approximation.
In Figure \ref{Ex1.error} the absolute error between the present method and exact solution for $N=10$ is plotted.
Figure \ref{Ex1.ai} represents the coefficients of the Bernoulli polynomials obtained by the present method for some various $N$ of the Bessel equation.
This figure shows that the method has an appropriate convergence rate.\\\\
\textbf{Example 2.}
In this example we consider the following Lane-Emden type equation \cite{Ramos.Laneemden,Chowdhury.Hashim.Laneemden,Bat.Hashim.Laneemden,parand.hermit.Laneemden}
\begin{eqnarray}\label{example2.main}
x~\frac{d^2u(x)}{dx^2}+8\frac{du(x)}{dx}+x^2~u(x)=x^6-x^5+44x^3-30x^2~,
\end{eqnarray}
with initial conditions
\begin{eqnarray}\label{initial.example2}
u(0)=0~,~~~~\frac{du(x)}{dx}|_{x=0}=0~,
\end{eqnarray}
The exact solution of this problem is
\begin{eqnarray}\nonumber
u(x)=x^4-x^3~.
\end{eqnarray}
Now, we approximate $d^2u(x)/dx^2$ by the Bernoulli polynomials as
\begin{eqnarray}\nonumber
\frac{d^2u(x)}{dx^2}=A^{T}\mathbf{B}(x)~,
\end{eqnarray}
Also, by using the initial conditions Eq. (\ref{initial.example2}) and the operation matrix of integration Eq. (\ref{oper.integral}), one has:
\begin{eqnarray}\nonumber
&&\frac{du(x)}{dx}=A^{T}\mathcal{I}\mathbf{B}(x)~,\\\label{exe2.y}
&&u(x)=A^{T}\mathcal{I}^{2}\mathbf{B}(x)~,
\end{eqnarray}
The functions $x,~ x^2,~x^6-x^5+44x^3-30x^2$ can be expressed as
\begin{eqnarray}\nonumber
&&x=E^{T}\mathbf{B}(x)~,\\\nonumber
&&x^2=F^{T}\mathbf{B}(x)~,\\\nonumber
&&x^6-x^5+44x^3-30x^2=V^{T}\mathbf{B}(x)~,
\end{eqnarray}
where
\begin{eqnarray}\nonumber
&&E=[1/2,~1,~0,...,~0]^{T}~,\\\nonumber
&&F=[1/3,~1,~1,~0,~0..,~0]^{T}~,\\\nonumber
&&V=[41/42,~14,~73/2,~137/3,~5/2,~2,~1,~0,~0,~...,~0]^{T}.
\end{eqnarray}
Therefore, equation (\ref{example2.main}) can be written as:
\begin{eqnarray}\nonumber
E^{T}\mathbf{B}(x)\mathbf{B}^{T}(x)A+8A^{T}\mathcal{I}\mathbf{B}(x)+F^{T}\mathbf{B}(x)\mathbf{B}^{T}(x)(\mathcal{I}^{2})^{T}A-V^{T}\mathbf{B}(x)=0~,
\end{eqnarray}
or
\begin{eqnarray}\nonumber
\Bigg(E^{T}\tilde{\mathbf{A}}+8A^{T}\mathcal{I}+F^{T}\tilde{\mathbf{Z}}-V^{T}\Bigg)\mathbf{B}(x)=0~,
\end{eqnarray}
where
\begin{eqnarray}\nonumber
&&\mathbf{B}(x)\mathbf{B}^{T}(x)A=\tilde{\mathbf{A}}\mathbf{B}(x)~,\\\nonumber
&&Z=(\mathcal{I}^{2})^{T}A~,\\\nonumber
&&\mathbf{B}(x)\mathbf{B}^{T}(x)Z=\tilde{\mathbf{Z}}\mathbf{B}(x)~,
\end{eqnarray}
By solving above equation for $N=4$, we obtain $A^T=\big(1,6,12\big)$.
Therefore
\begin{eqnarray}\nonumber
&&\frac{d^2u(x)}{dx^2}=A^{T}\mathbf{B}(x)=\big(1,6,12\big)
\left(
 \begin{array}{c}
1\\
x-1/2\\
x^2-x+1/6\\
\end{array}
\right)=12x^2-6x~,\\
&&\frac{du(x)}{dx}=\int_{0}^{x}\frac{d^2u(t)}{dt^2}dt=4x^3-3x^2,\\\nonumber
&&u(x)=\int_{0}^{x}\frac{du(t)}{dt}dt=x^4-x^3,\nonumber
\end{eqnarray}
which is the exact solution of the problem.\\\\
\textbf{Example 3.}
In this example we consider the following nonlinear Riccati equation \cite{Tawil.Bahnasawi.AbdelNaby.Riccati,Abbasbandy.Riccati,Abbasbandy.Riccati.2,Abbasbandy.Riccati.3,Geng.Lin.Cui.Riccati,FazhanGeng.Riccati}
\begin{eqnarray}\label{example3.main}
\frac{du(x)}{dx}=2u(x)-u^{2}(x)+1~,
\end{eqnarray}
with initial condition
\begin{eqnarray}\label{initial.example3}
u(0)=0~,
\end{eqnarray}
The exact solution of this problem is
\begin{eqnarray}\nonumber
u(x)=1+\sqrt{2}\tanh\Bigg(\sqrt{2}x+\frac{1}{2}\ln(\frac{\sqrt{2}-1}{\sqrt{2}+1})\Bigg)~.
\end{eqnarray}
Now, we approximate $du(x)/dx$ by the Bernoulli polynomials as
\begin{eqnarray}\nonumber
\frac{du(x)}{dx}=A^{T}\mathbf{B}(x)~,
\end{eqnarray}
Also, by using the initial conditions Eq. (\ref{initial.example3}) and the operation matrix of integration Eq. (\ref{oper.integral}), we have:
\begin{eqnarray}\label{exe3.y}
&&u(x)=A^{T}\mathcal{I}\mathbf{B}(x)~,\\
&&u^{2}(x)=A^{T}\mathcal{I}\mathbf{B}(x)\mathbf{B}(x)^{T}\mathcal{I}^{T}A=A^{T}\mathcal{I}\tilde{\mathbf{Z}}\mathbf{B}(x)~,
\end{eqnarray}
where
\begin{eqnarray}\nonumber
&&Z=\mathcal{I}^{T}A~,\\\nonumber
&&\mathbf{B}(x)\mathbf{B}(x)^{T}Z=\tilde{\mathbf{Z}}\mathbf{B}(x)~.
\end{eqnarray}
Therefore, equation (\ref{example3.main}) can be rewritten as:
\begin{eqnarray}\nonumber
-A^{T}+2A^{T}\mathcal{I}-A^{T}\mathcal{I}\tilde{\mathbf{Z}}+E^{T}=0~,
\end{eqnarray}
where
$1=E^{T}~\mathbf{B}(x)$.\\
As in a typical Galerkin method \cite{Gottlieb.Hussaini} we generate $N + 1$ equations, therefore by solving above equations the unknown vector $A$ is achieved and the unknown
$u(x)$ can be calculated by using Eq. (\ref{exe3.y}).\\
In Table \ref{Tab2}, a comparison is made between the
approximate values using the introduced technique
together with the exact solution for some various $N$.
Also, the RMS errors for some various $N$ are shown in Table \ref{Tab2}.
From Table \ref{Tab2}, it can be seen that
a few term of Bernoulli polynomials is sufficient to achieve a good approximation.
In Figure \ref{Ex3.error} the absolute error between our obtained approximate solutions and exact solution for $N=10$ is plotted.\\\\
\textbf{Example 4.}
Finally, we consider the following nonlinear Riccati equation \cite{Geng.Lin.Cui.Riccati,FazhanGeng.Riccati}
\begin{eqnarray}\label{example4.main}
\frac{du(x)}{dx}=1+x^2-u^{2}(x)~,
\end{eqnarray}
with initial condition
\begin{eqnarray}\label{initial.example4}
u(0)=1~,
\end{eqnarray}
and the exact solution of this problem is
\begin{eqnarray}\nonumber
u(x)=x+\frac{e^{-x^2}}{1+\int_{0}^{x}e^{-t^2}~dt}~.
\end{eqnarray}
Now, we approximate $du(x)/dx$ by the Bernoulli polynomials as
\begin{eqnarray}\nonumber
\frac{du(x)}{dx}=A^{T}\mathbf{B}(x)~,
\end{eqnarray}
Also, by using the initial conditions Eq. (\ref{initial.example3}) and the operation matrix of integration Eq. (\ref{oper.integral}), one has:
\begin{eqnarray}\label{exe4.y}
&&u(x)=A^{T}\mathcal{I}\mathbf{B}(x)+E^{T}\mathbf{B}(x)~,\\\nonumber
&&u^{2}(x)=\Bigg(A^{T}\mathcal{I}\tilde{\mathbf{Z}}+E^{T}\tilde{\mathbf{E}}+E^{T}\tilde{\mathbf{Z}}+A^{T}\mathcal{I}\tilde{\mathbf{E}}\Bigg)\mathbf{B}(x)~,
\end{eqnarray}
where
\begin{eqnarray}\nonumber
&&Z=\mathcal{I}^{T}A~,\\\nonumber
&&\mathbf{B}(x)\mathbf{B}(x)^{T}Z=\tilde{\mathbf{Z}}\mathbf{B}(x)~,\\\nonumber
&&\mathbf{B}(x)\mathbf{B}(x)^{T}E=\tilde{\mathbf{E}}\mathbf{B}(x)~.
\end{eqnarray}
The function $1+x^2$ is given as
\begin{eqnarray}\nonumber
&&1+x^2=F^{T}\mathbf{B}(x)~,
\end{eqnarray}
where
\begin{eqnarray}\nonumber
V=[4/3,~1,~1,~0,~0,~...,~0]^{T}.
\end{eqnarray}
Therefore, Eq. (\ref{example4.main}) can be rewritten as:
\begin{eqnarray}\nonumber
A^{T}-F^{T}+A^{T}\mathcal{I}\tilde{\mathbf{Z}}+E^{T}\tilde{\mathbf{E}}+E^{T}\tilde{\mathbf{Z}}+A^{T}\mathcal{I}\tilde{\mathbf{E}}=0~,
\end{eqnarray}
As in a typical Galerkin method \cite{Gottlieb.Hussaini} we generate $N + 1$ equations, therefore by solving above equations the unknown vector $A$ is achieved and the unknown
$u(x)$ can be calculated from Eq. (\ref{exe4.y}).\\
In Table \ref{Tab3}, a comparison is made between the
approximate values of the applied scheme
together with the exact solution for some various $N$.
Also, the RMS errors for some various $N$ are shown in Table \ref{Tab3}.
From Table \ref{Tab3}, it can be seen that
a few term of Bernoulli polynomials is sufficient to achieve a good approximation.
In Figure \ref{Ex4.error} the absolute error between our approximate results and the exact solution for $N=10$ is plotted.
Figure \ref{Ex4.ai} represents the coefficients of the Bernoulli polynomials obtained by the present method for some various $N$ of the Bessel equation.
This figure shows that the method has an appropriate convergence rate.\\\\
\section{Conclusions}\label{Conclusions}
Nonlinear differential and integral equations have a very important place in physics, mathematics
and engineering. Since this equations are usually difficult to solve analytically, it is required to obtain their approximate
solution. For this reason, the present method has been proposed to approximate the solutions of these equations using the
Bernoulli polynomials. In this paper, the Bernoulli polynomials operational matrices of integration, differentiation and product are derived. A
general procedure of forming these matrices are given.
The method is general, easy to implement,
and yields very accurate results. Moreover, only
a few number of basis yields a satisfactory
result.
\bibliographystyle{elsart}
\bibliography{Bernoulli}

\begin{thebibliography}{10}
\expandafter\ifx\csname url\endcsname\relax
  \def\url#1{\texttt{#1}}\fi
\expandafter\ifx\csname urlprefix\endcsname\relax\def\urlprefix{URL }\fi

\bibitem{Yildirim..Kaplan.HPM}
A.~Yildirim, S.~Sezer, Y.~Kaplan, Analytical approach to {B}oussinesq equation
  with space and time-fractional derivatives, Int. J. Numer. Meth. Fluids 66
  (2011) 1315--1324.

\bibitem{He.wave.HPM}
J.~H. He, Application of homotopy perturbation method to nonlinear wave
  equations, Chaos, Solitons and Fractals 26 (2005) 695--700.

\bibitem{Kazem.amani.RBF.EABE}
S.~Kazem, J.~Rad, K.~Parand, Radial basis functions methods for solving
  {F}okker-{P}lanck equation, Eng. Anal. Bound. Elem. (2011), in press.

\bibitem{parand.JCP}
K.~Parand, M.~Shahini, M.~Dehghan, Rational {L}egendre pseudospectral approach
  for solving nonlinear differential equations of {L}ane-{E}mden type, J.
  Comput. Phys. 228 (2009) 8830--8840.

\bibitem{Abbasbandy.Shivanian.HAM.INT}
S.~Abbasbandy, E.~Shivanian, A new analytical technique to solve {F}redholm's
  integral equations, Numer. Algorithms 56 (2011) 27--43.

\bibitem{Liao.HAM}
S.~J. Liao, Series solution of nonlinear eigenvalue problems by means of the
  homotopy analysis method, Nonlinear Analysis: Real World Applications 10
  (2009) 2455--2470.

\bibitem{Parand.Razzaghi.scripta}
K.~Parand, M.~Razzaghi, Rational {L}egendre approximation for solving some
  physical problems on semi-infinite intervals, Phys. Scr. 69 (2004) 353--357.

\bibitem{Parand.Razzaghi.IJCM}
K.~Parand, M.~Razzaghi, Rational {C}hebyshev tau method for solving
  higher-order ordinary differential equations, Int. J. Comput. Math. 81 (2004)
  73--80.

\bibitem{Yousefi.Dehghan.IJCM}
S.~Yousefi, M.~Dehghan, The use of {H}e's variational iteration method for
  solving variational problems, Int. J. Comput. Math. 87 (2010) 1299--1314.

\bibitem{wazwaz.VIM}
A.~Wazwaz, A reliable treatment of singular {E}mden-{F}owler initial value
  problems and boundary value problems, Appl. Math. Comput. 217 (2011)
  10387--10395.

\bibitem{Saadatmandi.Dehghan.Leg.Tau}
A.~Saadatmandi, M.~Dehghan, A new operational matrix for solving
  fractional-order differential equations, Comput. Math. Appl. 59 (2010)
  1326--1336.

\bibitem{Saadatmandi.Dehghan.Leg.Tau.2}
A.~Saadatmandi, M.~Dehghan, A tau approach for solution of the space fractional
  diffusion equation, Comput. Math. Appl. (2010) in press.

\bibitem{Gu.Jiang.Harr}
J.~S. Gu, W.~S. Jiang, The haar wavelets operational matrix of integration, I.
  J. Syst. Sci. 27 (1996) 623--628.

\bibitem{Horng.Chou}
I.~Horng, J.~Chou, Shifted {C}hebyshev direct method for solving variational
  problems, I. J. Syst. Sci. 16 (1985) 855--861.

\bibitem{Razzaghi.yousefi.wavelet}
M.~Razzaghi, S.~Yousefi, The legendre {W}avelets operational matrix of
  integration, I. J. Syst. Sci. 32 (2001) 495--502.

\bibitem{Razzaghi.yousefi.sin.cos}
M.~Razzaghi, S.~Yousefi, Sine-{C}osine {W}avelets operational matrix of
  integration and its applications in the calculus of variations, I. J. Syst.
  Sci. 33 (2002) 805--810.

\bibitem{S.A.Yousefi.Behroozifar.2010}
S.~A. Yousefi, M.~Behroozifar, Operational matrices of bernstein polynomials
  and their applications, I. J. Sys. Sci. 41 (2010) 709--716.

\bibitem{Norlund}
N.~E. Norlund, Vorlesungen uber Differenzenrechnung, Springer-Verlag, New York,
  1954.

\bibitem{Vandiver}
H.~Vandiver, Certain congruences involving the {B}ernoulli numbers, Duke
  Mathematical Journal 5 (1939) 548--551.

\bibitem{Cheon.Bernoulli}
G.~Cheon, A note on the {B}ernoulli and {E}uler polynomials, Appl. Math. Lett.
  16 (2003) 365--368.

\bibitem{Burak.Simsek.Bernoul}
B.~Kurt, Y.~Simsek, Notes on generalization of the {B}ernoulli type
  polynomials, Appl. Math. Comput. (2011) doi:10.1016/j.amc.2011.03.086.

\bibitem{Lu.Bernoulli}
D.~Lu, Some properties of {B}ernoulli polynomials and their generalizations,
  Appl. Math. Lett. 24 (2011) 746--751.

\bibitem{Agoh.Dilcher.Bernoulli}
T.~Agoh, K.~Dilcher, Integrals of products of {B}ernoulli polynomials, J. Math.
  Anal. Appl. 381 (2011) 10--16.

\bibitem{Buric.Elezovic.Bernoulli}
T.~Buric, N.~Elezovic, Bernoulli polynomials and asymptotic expansions of the
  quotient of gamma functions, J. Comput. Appl. Math. 235 (2011) 3315--3331.

\bibitem{Natalini.Bernoulli}
P.~Natalini, A.~Bernardini, A generalization of the {B}ernoulli polynomials, J.
  Appl. Math. 2003 (2003) 155--163.

\bibitem{Kurt.Bernoulli}
B.~Kurt, A further generalization of the {B}ernoulli polynomials and on the
  $2d$-{B}ernoulli polynomials ${B}_{n}^{2}(x, y)$, Appl. Math. Sci. 4 (2010)
  2315--2322.

\bibitem{Gottlieb.Hussaini}
D.~Gottlieb, M.~Hussaini, S.~Orszg, Theory and applications of spectral methods
  in spectral methods for partial differential equations, SIAM, Philadelphia,
  1984.

\bibitem{Costabile}
F.~Costabile, Expansions of real functions in {B}ernoulli polynomials and
  applications, Conf. Sem. Mat.Univ. Bari, N. 273 (1999) 1--13.

\bibitem{ONeil}
P.~O'Neil, Advanced Engineering Mathematics, Belmont, California, 1987.

\bibitem{Yousefi.Behroozifar.IJSS}
S.~Yousefi, M.~Behroozifar, Operational matrices of {B}ernstein polynomials and
  their applications, I. J. Sys. Sci. 41 (2010) 709--716.

\bibitem{Ramos.Laneemden}
J.~Ramos, Linearization techniques for singular initial-value problems of
  ordinary differential equations, Appl. Math. Comput. 161 (2005) 525--542.

\bibitem{Chowdhury.Hashim.Laneemden}
M.~Chowdhury, I.~Hashim, Solutions of {E}mden-{F}owler equations by homotopy
  perturbation method, Nonlinear Anal. Real World Appl. 10 (2009) 104--115.

\bibitem{Bat.Hashim.Laneemden}
A.~Bataineh, M.~Noorani, I.~Hashim, Homotopy analysis method for singular
  {IVP}s of {E}mden-{F}owler type, Commun. Nonlinear Sci. Numer. Simul. 14
  (2009) 1121--1131.

\bibitem{parand.hermit.Laneemden}
K.~Parand, M.~Dehghan, A.~Rezaei, S.~Ghaderi, An approximation algorithm for
  the solution of the nonlinear {L}ane-{E}mden type equations arising in
  astrophysics using {H}ermite functions collocation method, Comput. Phys.
  Commun. 181 (2010) 1096--1108.

\bibitem{Tawil.Bahnasawi.AbdelNaby.Riccati}
M.~El-Tawil, A.~Bahnasawi, A.~Abdel-Naby, Solving {R}iccati differential
  equation using {A}domian's decomposition method, Appl. Math. Comput. 157
  (2004) 503--514.

\bibitem{Abbasbandy.Riccati}
S.~Abbasbandy, Homotopy perturbation method for quadratic {R}iccati
  differential equation and comparison with {A}domian's decomposition method,
  Appl. Math. Comput. 172 (2006) 485--490.

\bibitem{Abbasbandy.Riccati.2}
S.~Abbasbandy, A new application of {H}e's variational iteration method for
  quadratic {R}iccati differential equation by using {A}domian's polynomials,
  J. Comput. Appl. Math. 207 (2007) 59--63.

\bibitem{Abbasbandy.Riccati.3}
S.~Abbasbandy, Iterated {H}e's homotopy perturbation method for quadratic
  {R}iccati differential equation, Appl. Math. Comput. 175 (2006) 581--589.

\bibitem{Geng.Lin.Cui.Riccati}
F.~Geng, Y.~Lin, M.~Cui, A piecewise variational iteration method for {R}iccati
  differential equations, Comput. Math. Appl. 58 (2009) 2518--2522.

\bibitem{FazhanGeng.Riccati}
F.~Geng, A modified variational iteration method for solving {R}iccati
  differential equations, Comput. Math. Appl. 60 (2010) 1868--1872.

\end{thebibliography}

\clearpage

\clearpage
\begin{table}\tiny
\caption{A comparison between the approximate values of the present approach together with the exact solution and the RMS errors for some various $N$ in Example 1.}\vline
\begin{tabular*}{\columnwidth}{@{\extracolsep{\fill}}*{12}{c}}
\hline
~~$x$ &\vline & Exact & $N=2$& $N=4$ & $N=6$& $N=8$ & $N=10$&\vline\\
\hline
~~~$0.0$&\vline & $1.0000000000$ & $1.0014754886$ & $0.9999950817$ & $1.0000000075$ & $1.0000000000$ & $1.0000000000$&\vline\\
~~~$0.2$&\vline & $0.9900249722$ & $0.9895787397$ & $0.9900259566$ & $0.9900249747$ & $0.9900249722$ & $0.9900249722$&\vline\\
~~~$0.4$&\vline & $0.9603982267$ & $0.9598567615$ & $0.9603965639$ & $0.9603982244$ & $0.9603982267$ & $0.9603982267$&\vline\\
~~~$0.6$&\vline & $0.9120048635$ & $0.9123095542$ & $0.9120062526$ & $0.9120048657$ & $0.9120048635$ & $0.9120048635$&\vline\\
~~~$0.8$&\vline & $0.8462873527$ & $0.8469371176$ & $0.8462868378$ & $0.8462873503$ & $0.8462873527$ & $0.8462873527$&\vline\\
~~~$1.0$&\vline & $0.7651976866$ & $0.7637394519$ & $0.7652025998$ & $0.7651976790$ & $0.7651976866$ & $0.7651976866$&\vline\\
\hline
~~~RMS &\vline &       $-$       & $7.9298 \times 10^{-4}$ & $2.5651 \times 10^{-6}$ & $3.7742 \times 10^{-9}$ & $3.1762 \times 10^{-12}$ & $1.901 \times 10^{-15}$&\vline\\
\hline
\end{tabular*}
\label{Tab1}
\end{table}

\begin{table}\tiny
\caption{A comparison between the approximate values of the present approach together with the exact solution and the RMS errors for some various $N$ in Example 3.}\vline
\begin{tabular*}{\columnwidth}{@{\extracolsep{\fill}}*{12}{c}}
\hline
~~$x$ &\vline & Exact & $N=4$& $N=6$ & $N=8$& $N=10$ & $N=14$&\vline\\
\hline
~~~$0.0$&\vline & $0.0000000000$ & $0.0025265013$ & $0.0001128571$ & $0.0000044867$ & $0.0000000000$ & $0.0000000000$&\vline\\
~~~$0.2$&\vline & $0.2419767996$ & $0.2412798535$ & $0.2419431649$ & $0.2419774162$ & $0.2419768508$ & $0.2419767996$&\vline\\
~~~$0.4$&\vline & $0.5678121663$ & $0.5687080628$ & $0.5678457562$ & $0.5678129586$ & $0.5678121631$ & $0.5678121662$&\vline\\
~~~$0.6$&\vline & $0.9535662164$ & $0.9529191921$ & $0.9535337210$ & $0.9535648227$ & $0.9535661622$ & $0.9535662164$&\vline\\
~~~$0.8$&\vline & $1.3463636554$ & $1.3465915432$ & $1.3464020534$ & $1.3463645265$ & $1.3463636335$ & $1.3463636553$&\vline\\
~~~$1.0$&\vline & $1.6894983916$ & $1.6869736561$ & $1.6896112463$ & $1.6894939049$ & $1.6894985427$ & $1.6894983916$&\vline\\
\hline
~~~RMS &\vline &       $-$       & $1.3164 \times 10^{-3}$ & $5.7405 \times 10^{-5}$ & $2.2057 \times 10^{-6}$ & $7.7612 \times 10^{-8}$ & $6.9755 \times 10^{-11}$&\vline\\
\hline
\end{tabular*}
\label{Tab2}
\end{table}

\begin{table}\tiny
\caption{A comparison between the approximate values of the present approach together with the exact solution and the RMS errors for some various $N$ in Example 4.}\vline
\begin{tabular*}{\columnwidth}{@{\extracolsep{\fill}}*{12}{c}}
\hline
~~$x$ &\vline & Exact & $N=2$& $N=4$ & $N=6$& $N=8$ & $N=10$&\vline\\
\hline
~~~$0.0$&\vline & $1.0000000000$ & $1.0058680884$ & $1.0000458868$ & $1.0000021902$ & $0.9999999813$ & $1.0000000000$&\vline\\
~~~$0.2$&\vline & $1.0024198255$ & $0.9999551044$ & $1.0024189537$ & $1.0024205376$ & $1.0024198267$ & $1.0024198255$&\vline\\
~~~$0.4$&\vline & $1.0176508789$ & $1.0168990769$ & $1.0176602221$ & $1.0176502287$ & $1.0176508744$ & $1.0176508789$&\vline\\
~~~$0.6$&\vline & $1.0544668099$ & $1.0567000057$ & $1.0544519892$ & $1.0544674572$ & $1.0544668142$ & $1.0544668099$&\vline\\
~~~$0.8$&\vline & $1.1180925454$ & $1.1193578910$ & $1.1181039457$ & $1.1180918378$ & $1.1180925448$ & $1.1180925453$&\vline\\
~~~$1.0$&\vline & $1.2105990147$ & $1.2048727327$ & $1.2105531747$ & $1.2105968244$ & $1.2105990333$ & $1.2105990151$&\vline\\
\hline
~~~RMS &\vline &       $-$       & $3.0751 \times 10^{-3}$ & $2.3575 \times 10^{-5}$ & $1.0961 \times 10^{-6}$ & $8.8580 \times 10^{-9}$ & $2.0730 \times 10^{-10}$&\vline\\
\hline
\end{tabular*}
\label{Tab3}
\end{table}
\end{document}